\documentclass[12pt,oneside]{amsart}

\usepackage{setspace}
\usepackage{amssymb}   
\usepackage{amsmath}
\usepackage{mathrsfs}
\usepackage{stmaryrd}
\usepackage{amsthm}
\usepackage{newlfont}
\usepackage{amscd}
\usepackage{mathtools}
\usepackage{bm}
\usepackage{tikz-cd}
\usepackage{graphicx}
\usepackage{hyperref}
\usepackage{enumitem}
\usepackage[all]{xy}
\usepackage{textcomp}
\usepackage{amsbsy}

\newtheorem{thm}{Theorem}[section]
\newtheorem{thm-defn}[thm]{Theorem/Definition}
\newtheorem{lem}[thm]{Lemma}
\newtheorem{prop}[thm]{Proposition}

\newtheorem{conj}[thm]{Conjecture}

\theoremstyle{definition}
\newtheorem{defn}[thm]{Definition}

\theoremstyle{remark}
\newtheorem{rem}[thm]{Remark}

\numberwithin{equation}{section}

\usepackage{relsize}
\usepackage[bbgreekl]{mathbbol}
\usepackage{amsfonts}

\DeclareSymbolFontAlphabet{\mathbb}{AMSb}
\DeclareSymbolFontAlphabet{\mathbbl}{bbold}
\newcommand{\Prism}{{\mathlarger{\mathbbl{\Delta}}}}

\begin{document}

\pagenumbering{arabic}

\title{On Fontaine's conjecture for torsion crystalline local systems}
\author{Yong Suk Moon}
\address{Beijing Institute of Mathematical Sciences and Applications, Beijing 101408, China}
\email{ysmoon@bimsa.cn}

\begin{abstract}
Let $\mathfrak{X}$ be a smooth connected $p$-adic formal scheme. Based on the prismatic description of crystalline local systems, we prove an analogue of Fontaine's conjecture for torsion crystalline local systems on the generic fiber of $\mathfrak{X}$. As an application, we show that the locus of crystalline local systems whose Hodge--Tate weights lie in a fixed interval cuts out a closed subscheme of the universal deformation ring.
\end{abstract}

\maketitle

\tableofcontents

\section{Introduction} \label{sec:intro}

Let $p$ be a prime. In \cite{fontaine-mazur}, Fontaine--Mazur made a deep conjecture that if a global $p$-adic Galois representation is potentially semi-stable at places dividing $p$ and unramified outside finitely many places, then it comes from algebraic geometry. One motivating aspect of the conjecture originated from studying the Galois deformation space satisfying certain $p$-adic Hodge-theoretic conditions. In particular, Fontaine (\cite{fontaine-deforming}) predicted that there exists a quotient of the universal deformation ring parametrizing crystalline representations. 

Fontaine's prediction can be formulated in terms of torsion representations as follows. Let $K$ be a complete discrete valuation field of mixed characteristic $(0, p)$ with perfect residue field $k$, and denote $G_K \coloneqq \mathrm{Gal}(\overline{K}/K)$.

\begin{conj}[cf.~{\cite{fontaine-deforming}}] \label{conj:Fontaine}
Fix a postivie integer $r$, and let $T$ be a finite free $\mathbf{Z}_p$-representation of $G_K$. Suppose for each $n \geq 1$, $T/p^n T$ is torsion crystalline with Hodge--Tate weights in $[0, r]$, i.e. there exist $G_K$-stable $\mathbf{Z}_p$-lattices $T_2^{(n)} \subset T_1^{(n)}$ of a crystalline representation with Hodge--Tate weights in $[0, r]$ such that $T/p^nT \cong T_1^{(n)} / T_2^{(n)}$ as $\mathbf{Z}_p[G_K]$-modules. Then $T[p^{-1}]$ is crystalline with Hodge--Tate weights in $[0, r]$. 	
\end{conj}

Conjecture \ref{conj:Fontaine} is proved by Liu (\cite{liu-fontaineconjecture}). Furthermore, Kisin studied the geometric structure of the deformation space associated with more refined $p$-adic Hodge-theoretic conditions (\cite{kisin-semistabledeformation}), and used it in the proof of Fontaine--Mazur conjecture for many cases of two dimensional representations (\cite{kisin-fontaine-mazur}).

For \'etale isogeny $\mathbf{Z}_p$-local systems (i.e. elements in the isogeny category of $\mathbf{Z}_p$-local systems as in \cite[Def.~1.4.1]{kedlaya-liu-relative-padichodge}) on a geometrically connected algebraic variety over a number field, the relative version of Fontaine--Mazur conjecture is formulated in \cite{liu-zhu-rigidity}. This turns out to be closely related to the expected compatibility between the classical Riemann--Hilbert correspondence and its $p$-adic analogue, and such a compatibility is proved in the case of Shimura varieties (\cite{diao-lan-liu-zhu-logRH}). Motivated by these results, we study in this paper a statement analogous to Conjecture \ref{conj:Fontaine} for geometric families of Galois representations.   

Let $\mathfrak{X}$ be a smooth connected $p$-adic formal scheme over $\mathcal{O}_K$ and $X$ be the generic fiber of $\mathfrak{X}$. Choose a geometric base point $\overline{x} \in X$, and let $G_X \coloneqq \pi_1(X, \overline{x})$ be the profinite fundamental group as in \cite[\S 3]{scholze-p-adic-hodge}. Note that by \cite[Prop.~3.5, 8.2]{scholze-p-adic-hodge}, the category of \'etale $\mathbf{Z}_p$-local systems on $X$ is naturally equivalent to the category $\mathrm{Rep}_{\mathbf{Z}_p}^{\mathrm{free}}(G_X)$ of continuous finite free $\mathbf{Z}_p$-representations of $G_X$. The notion of crystalline \'etale $\mathbf{Z}_p$-local systems on $X$ are studied in \cite{du-liu-moon-shimizu-completed-prismatic-F-crystal-loc-system} and \cite{guo-reinecke-prism-F-crys} (see also \cite{brinon-relative} for small affine cases).

\begin{thm}[Thm.~\ref{thm:main}] \label{thm:intro}
Fix a positive integer $r$, and let $T \in \mathrm{Rep}_{\mathbf{Z}_p}^{\mathrm{free}}(G_X)$. Suppose for each $n \geq 1$, $T/p^n T$ is torsion crystalline with Hodge--Tate weights in $[0, r]$. Then $T\otimes_{\mathbf{Z}_p} \mathbf{Q}_p$ is a crystalline local system on $X$ with Hodge--Tate weights in $[0, r]$.
\end{thm}

\noindent We remark that unlike the classical situation $\mathfrak{X} = \mathrm{Spf} (\mathcal{O}_K)$, the corresponding statement for Barsotti--Tate representations does not hold in general (\cite{moon-relativeRaynaud}). So the relative situation turns out to be more subtle. 

As an application, we study the corresponding crystalline deformation ring. For a fixed absolutely irreducible $\mathbf{F}_p$-representation $V_0$ of $G_X$, there is a universal deformation ring parametrizing deformations of $V_0$ (\cite{smit-lenstra}).

\begin{thm}[Thm.~\ref{thm:cryst-deform-ring}] \label{thm:cryst-deform-ring-intro}
Fix an integer $r  > 0$. The locus of crystalline local systems on $X$ with Hodge--Tate weights in $[0, r]$ cuts out a closed subscheme of the universal deformation scheme. 	
\end{thm}

\begin{rem}
Let $\mathfrak{X}$ be a connected \textit{semistable} $p$-adic formal scheme over $\mathcal{O}_K$ and $X$ be the generic fiber of $\mathfrak{X}$. Then by the purity of semistable local systems recently proved in \cite{du-liu-moon-shimizu-log-prism-F-cryst-purity}, the same arguments as in this paper would imply that the statements analogous to Theorem~\ref{thm:intro} and Theorem~\ref{thm:cryst-deform-ring-intro} for \textit{semistable} local systems on $X$ also hold. 
\end{rem}

\noindent \textbf{Strategy of the proof.} We briefly explain the main ideas in the proof of Theorem \ref{thm:intro}. By the local description of crystalline local systems (\cite[\S A]{du-liu-moon-shimizu-completed-prismatic-F-crystal-loc-system}), we can reduce to the case when $\mathfrak{X} = \mathrm{Spf} R$ where $R$ is the $p$-completion of an \'etale algebra over $\mathcal{O}_K\langle t_1^{\pm 1}, \ldots, t_d^{\pm 1}\rangle$ for some $d \geq 0$. Let $\pi \in K$ be a uniformizer. Using the purity result for crystalline representations (\cite{moon-purity}, see also \cite[\S 5.4]{tsuji-cryst-shvs}), we further reduce to showing the analogous statement for representations of $G_L \coloneqq \mathrm{Gal}(\overline{L}/L)$ where $L$ is the $p$-completion of $R_{(\pi)}$. The category of lattices of crystalline representations of $G_L$ is equivalent to the category of prismatic $F$-crystals on the absolute prismatic site of $\mathrm{Spf} \mathcal{O}_L$, or equivalently, that of Kisin modules with descent data (\cite[Cor.~4.36]{du-liu-moon-shimizu-completed-prismatic-F-crystal-loc-system}, see also \cite{guo-reinecke-prism-F-crys}). We obtain the corresponding Kisin module via a similar method as in \cite[\S 4.3]{liu-fontaineconjecture}. Then we construct the associated descent datum based on the equivalence between $\mathbf{Z}_p$-representations of $G_L$ and Laurent $F$-crystals at the derived level (\cite[\S 3]{bhatt-scholze-prismaticFcrystal}).\\ 

\noindent \textbf{Outline.} In Section 2, the purity result for crystalline local systems is explained. In Section 3, we review the theory of \'etale $\varphi$-modules and Galois representations. Section 4 discusses the prismatic description of crystalline $\mathbf{Z}_p$-local systems in terms of Kisin descent data. We prove Theorem~\ref{thm:intro} in Section 5, based on the results in previous sections. In Section 6, we study the crystalline deformation ring and show Theorem~\ref{thm:cryst-deform-ring-intro} as an application of Theorem~\ref{thm:intro}.\\   

\noindent \textbf{Notation and convention.} Let $p$ be a prime and $k$ be a perfect field of characteristic $p$. Let $K$ be a finite totally ramified extension of $K_0 \coloneqq W(k)[p^{-1}]$, and let $\mathcal{O}_K$ be its ring of integers. Choose a uniformizer $\pi \in K$ and let $E = E(u) \in W(k)[u]$ denote the monic minimal polynomial of $\pi$. Write $G_K \coloneqq \mathrm{Gal}(\overline{K}/K)$ for the absolute Galois group of $K$.

Let $\mathfrak{X}$ be a smooth connected $p$-adic formal scheme over $\mathcal{O}_K$ and $X$ be the generic fiber of $\mathfrak{X}$. Choose a geometric base point $\overline{x} \in X$, and write $G_X \coloneqq \pi_1(X, \overline{x})$ for the profinite fundamental group as in \cite[\S 3]{scholze-p-adic-hodge}. Let $\mathrm{Rep}_{\mathbf{Z}_p}^{\mathrm{free}}(G_X)$ be the category of continuous finite free $\mathbf{Z}_p$-representations of $G_X$. By a representation, we always mean a continuous representation. 

For any $\mathbf{Z}$-module $A$, we denote by $A^{\wedge}_p$ the classical $p$-completion of $A$. In our convention, the cyclotomic character has Hodge--Tate weight one.

\section*{Acknowledgements}

I would like to thank Mark Kisin and Tong Liu for their valuable comments on an earlier version of this article.

\section{Purity of crystalline representations} \label{sec:purity}

We explain the purity result for crystalline local systems proved in \cite{moon-purity}. First, note that our $p$-adic formal scheme $\mathfrak{X}$ can be covered Zariski locally by $\mathrm{Spf} R$, where $R$ is the $p$-completion of an \'etale algebra over $\mathcal{O}_K\langle t_1^{\pm 1}, \ldots, t_d^{\pm 1}\rangle$ for some $d \geq 0$. For Theorem~\ref{thm:intro}, it suffices to consider this local case $\mathfrak{X} = \mathrm{Spf} R$ by \cite[Prop.~A.10]{du-liu-moon-shimizu-completed-prismatic-F-crystal-loc-system}. 

Write $G_{R[p^{-1}]}$ for the \'etale fundamental group $\pi_1^{\text{\'et}}(\mathrm{Spec} (R[p^{-1}]), \overline{x})$. By \cite[Ex.~1.6.6 ii)]{Huber-etale} and \cite[Rem.~1.4.4]{kedlaya-liu-relative-padichodge}, the category of $\mathbf{Z}_p$-local systems on the \'etale site of the adic space $\mathrm{Spa}(R[p^{-1}], R)$ is naturally equivalent to the category of finite free $\mathbf{Z}_p$-modules with continuous $G_{R[p^{-1}]}$-action. Note that $\mathrm{Spa}(R[p^{-1}], R)$ is the adic generic fiber of $\mathrm{Spf} R$. 

Let $R_0 \subset R$ be a $W(k)$-subalgebra which is the $p$-completion of an \'etale algebra over $W(k)\langle t_1^{\pm 1}, \ldots, t_d^{\pm 1}\rangle$ such that $R = R_0\otimes_{W(k)} \mathcal{O}_K$. Let $L_0$ be the $p$-completion of $R_{0, (p)}$. $L_0$ is a complete discrete valuation ring whose residue field has a finite $p$-basis given by $\{t_1, \ldots, t_d\}$. We equip $L_0$ with the Frobenius endomorphism such that $\varphi(t_i) = t_i^p$, and let $\mathcal{O}_{L_0}$ be its ring of integers. Let $L = L_0\otimes_{W(k)} \mathcal{O}_K$ and denote by $\mathcal{O}_L$ its rings of integers. Write $\overline{L}$ for an algebraic closure of $L$. Let $V$ be a finite $\mathbf{Q}_p$-representation of $G_{R[p^{-1}]}$. By the change of paths for \'etale fundamental groups, we get a continuous group homomorphism $G_L\coloneqq \mathrm{Gal}(\overline{L}/L) \rightarrow G_{R[p^{-1}]}$. So $V$ can also be considered as a representation of $G_L$. In \cite{moon-purity} (see also \cite[Thm.~5.4.8]{tsuji-cryst-shvs}), the following purity result is proved.

\begin{thm}[cf. {\cite[Thm.~1.2]{moon-purity}}] \label{thm:purity}
Let $r > 0$ be an integer. If $V|_{G_L}$ is a crystalline representation of $G_L$ with Hodge--Tate weights in $[0, r]$, then $V$ is a crystalline representation of $G_{R[p^{-1}]}$ with Hodge--Tate weights in $[0, r]$.	
\end{thm}
 
\noindent Hence, to prove Theorem~\ref{thm:intro}, it suffices to show the corresponding statement for representations of $G_L$. 

The formalism of crystalline representations of $G_L$ is studied in \cite{brinon-crys-rep-imperfect-residue} using the the crystalline period ring $B_{\mathrm{cris}}$ defined in \cite[\S 2.3]{brinon-crys-rep-imperfect-residue}, which generalizes the classical case. For any finite $\mathbf{Q}_p$-representation $V$ of $G_L$, write 
\[
D_{\mathrm{cris}}(V) \coloneqq (V\otimes_{\mathbf{Q}_p} B_{\mathrm{cris}})^{G_L}
\]
which is a module over $L_0$. By the results proved in \cite[\S 3]{brinon-crys-rep-imperfect-residue}, we have $\dim_{L_0} D_{\mathrm{cris}}(V) \leq \dim_{\mathbf{Q}_p} V$, and the equality holds if and only if $V$ is crystalline. The following Lemmas also hold as in the classical case and will be used in Section~\ref{sec:cryst-deform-ring}. 

\begin{lem} \label{lem:exact-seq-reps}
Let $V$ be a finite $\mathbf{Q}_p$-representation of $G_L$, and let $W \subset V$ be a sub-representation. If $V$ is crystalline, then $W$ is crystalline. 	
\end{lem}

\begin{proof}
Consider the exact sequence of representations
\[
0 \rightarrow W \rightarrow V \rightarrow V/W \rightarrow 0.
\]
This induces an exact sequence of $L_0$-modules
\[
0 \rightarrow D_{\mathrm{cris}}(W) \rightarrow D_{\mathrm{cris}}(V) \rightarrow D_{\mathrm{cris}}(V/W).
\]
We have $\dim_{L_0} D_{\mathrm{cris}}(V) = \dim_{\mathbf{Q}_p}(V)$, $\dim_{L_0} D_{\mathrm{cris}}(W) \leq \dim_{\mathbf{Q}_p} W$, and $\dim_{L_0} D_{\mathrm{cris}}(V/W) \leq \dim_{\mathbf{Q}_p}(V/W)$. Since $\dim_{\mathbf{Q}_p}(V) = \dim_{\mathbf{Q}_p}(W)+\dim_{\mathbf{Q}_p}(V/W)$, we deduce that $\dim_{L_0} D_{\mathrm{cris}}(W) = \dim_{\mathbf{Q}_p} W$, so $W$ is crystalline.
\end{proof}

\begin{lem} \label{lem:base-change-coeff-reps}
Let $B$ be a finite $\mathbf{Q}_p$-algebra, and let $V_B$ be finite free $B$-module with a $B$-linear $G_L$-action which is continuous (considering $V_B$ as a $\mathbf{Q}_p$-representation). Let $B'$ be a finite $B$-algebra. Suppose $V_B$ is crystalline. Then the $B'$-linear representation $B'\otimes_B V_B$ of $G_L$ is crystalline, and
\[
D_{\mathrm{cris}}(B'\otimes_B V_B) = B'\otimes_B D_{\mathrm{cris}}(V_B)
\]	
\end{lem}

\begin{proof}
Since $V_B$ is crystalline, the natural map
\[
D_{\mathrm{cris}}(V_B)\otimes_{L_0} B_{\mathrm{cris}} \rightarrow V_B\otimes_{\mathbf{Q}_p} B_{\mathrm{cris}}
\]	
is an isomorphism. This induces an isomorphism
\[
(B'\otimes_B D_{\mathrm{cris}}(V_B))\otimes_{L_0} B_{\mathrm{cris}} \cong (B'\otimes_B V_B)\otimes_{\mathbf{Q}_p} B_{\mathrm{cris}}.
\]
In particular, we have $B'\otimes_B D_{\mathrm{cris}}(V_B) \subset D_{\mathrm{cris}}(B'\otimes_B V_B)$. On the other hand,
\[
\dim_{L_0}D_{\mathrm{cris}}(B'\otimes_B V_B) \leq \dim_{\mathbf{Q}_p} (B'\otimes_B V_B) = \dim_{L_0} (B'\otimes_B D_{\mathrm{cris}}(V_B))
\]
where the equality follows from the above isomorphism. Thus,
\[
B'\otimes_B D_{\mathrm{cris}}(V_B) = D_{\mathrm{cris}}(B'\otimes_B V_B)
\]
and $B'\otimes_B V_B$ is crystalline.
\end{proof}

\section{\'Etale $\varphi$-modules}

We briefly recall the theory of \'etale $\varphi$-modules and relations to Galois representations. Let $\mathfrak{S} \coloneqq \mathcal{O}_{L_0}[\![u]\!]$ equipped with the Frobenius extending that on $\mathcal{O}_{L_0}$ by $\varphi(u) = u^p$, and let $\mathcal{O}_{\mathcal{E}} \coloneqq \mathfrak{S}[u^{-1}]^{\wedge}_p$ equipped with the Frobenius extending $\varphi$ on $\mathfrak{S}$.   

\begin{defn} \label{def:etale-phi-mod}
An \textit{\'etale $(\varphi, \mathcal{O}_{\mathcal{E}})$-module} is a pair $(\mathcal{M}, \varphi_{\mathcal{M}})$ where $\mathcal{M}$ is a finitely generated $\mathcal{O}_{\mathcal{E}}$-module and $\varphi_{\mathcal{M}}\colon \mathcal{M} \rightarrow \mathcal{M}$ is a $\varphi$-semi-linear endomorphism such that $1\otimes\varphi_{\mathcal{M}}\colon \varphi^*\mathcal{M} = \mathcal{O}_{\mathcal{E}}\otimes_{\varphi, \mathcal{O}_{\mathcal{E}}}\mathcal{M} \rightarrow \mathcal{M}$ is an isomorphism. We say that an \'etale $(\varphi, \mathcal{O}_{\mathcal{E}})$-module is \textit{free} (resp. \textit{torsion}) if the underlying $\mathcal{O}_{\mathcal{E}}$-module $\mathcal{M}$ is free (resp. $p$-power torsion). 	

Let $\mathrm{Mod}_{\mathcal{O}_{\mathcal{E}}}$ denote the category of \'etale $(\varphi, \mathcal{O}_{\mathcal{E}})$-modules whose morphisms are $\mathcal{O}_{\mathcal{E}}$-linear maps compatible with Frobenii. Let $\mathrm{Mod}_{\mathcal{O}_{\mathcal{E}}}^{\mathrm{free}}$ and $\mathrm{Mod}_{\mathcal{O}_{\mathcal{E}}}^{\mathrm{tor}}$ denote respectively the full subcategories of free and torsion objects.   
\end{defn}

We use \'etale $(\varphi, \mathcal{O}_{\mathcal{E}})$-modules to study certain Galois representations. Let $\mathcal{O}_{\overline{L}}$ denote the ring of integers of $\overline{L}$. For integers $n \geq 0$, choose $\pi_n \in \overline{K}$ compatibly such that $\pi_0 = \pi$ and $\pi_{n+1}^p = \pi_n$, and let $\underline{\pi} = (\pi_n)_{n \geq 0} \in \mathcal{O}_{\overline{L}}^{\flat} = \varprojlim_{\varphi} \mathcal{O}_{\overline{L}}/p\mathcal{O}_{\overline{L}}$. For each $i = 1, \ldots, d$, choose $t_{i, n} \in \mathcal{O}_{\overline{L}}$ such that $t_{i, 0} = t_i$ and $t_{i, n+1}^p = t_{i, n}$, and write $\underline{t_i} = (t_{i, n})_{n \geq 0} \in \mathcal{O}_{\overline{L}}^{\flat}$. Let $\widetilde{L}_{\infty} \coloneqq \bigcup_{n \geq 0} L(\pi_n, t_{1, n}, \ldots, t_{d, n})$, and let $G_{\widetilde{L}_{\infty}} \coloneqq \mathrm{Gal}(\overline{L}/\widetilde{L}_{\infty})$ be the Galois subgroup of $G_L$. Consider the $W(k)$-linear map $\mathcal{O}_{L_0} \rightarrow W(\mathcal{O}_{\overline{L}}^{\flat})$ which maps $t_i$ to $[\underline{t_i}]$ and is compatible with Frobenius. This induces a $\varphi$-equivariant embedding $\mathfrak{S} \rightarrow W(\mathcal{O}_{\overline{L}}^{\flat})$ given by $u \mapsto [\underline{\pi}]$, which extends to $\mathcal{O}_{\mathcal{E}} \hookrightarrow W(\overline{L}^{\flat})$. Let $\widehat{\mathcal{O}}_{\mathcal{E}}^{\mathrm{ur}}$ be the $p$-completion of the ring of integers of the maximal unramified extension of $\mathcal{O}_{\mathcal{E}}[p^{-1}]$ inside $W(\overline{L}^{\flat})[p^{-1}]$. By \cite[Lem.~7.5, 7.6]{kim-groupscheme-relative}, we have $(\widehat{\mathcal{O}}_{\mathcal{E}}^{\mathrm{ur}})^{G_{\widetilde{L}_{\infty}}} = \mathcal{O}_{\mathcal{E}}$, and there exists a unique $G_{\widetilde{L}_{\infty}}$-equivariant ring endomorphism $\varphi$ on $\widehat{\mathcal{O}}_{\mathcal{E}}^{\mathrm{ur}}$ lifting the $p$-th power map on $\widehat{\mathcal{O}}_{\mathcal{E}}^{\mathrm{ur}}/(p)$ and lifting $\varphi$ on $\mathcal{O}_{\mathcal{E}}$. Furthermore, the map $\widehat{\mathcal{O}}_{\mathcal{E}}^{\mathrm{ur}} \hookrightarrow W(\overline{L}^{\flat})$ is compatible with $\varphi$.

Write $\mathrm{Rep}_{\mathbf{Z}_p}(G_{\widetilde{L}_{\infty}})$ for the category of finite $\mathbf{Z}_p$-modules with continuous $G_{\widetilde{L}_{\infty}}$-action, and let $\mathrm{Rep}_{\mathbf{Z}_p}^{\mathrm{free}}(G_{\widetilde{L}_{\infty}})$ and $\mathrm{Rep}_{\mathbf{Z}_p}^{\mathrm{tor}}(G_{\widetilde{L}_{\infty}})$ denote respectively the full subcategories of free and torsion objects. For $\mathcal{M} \in \mathrm{Mod}_{\mathcal{O}_{\mathcal{E}}}$ and $T \in \mathrm{Rep}_{\mathbf{Z}_p}(G_{\widetilde{L}_{\infty}})$, define
\[
T(\mathcal{M}) \coloneqq (\mathcal{M}\otimes_{\mathcal{O}_{\mathcal{E}}}\widehat{\mathcal{O}}_{\mathcal{E}}^{\mathrm{ur}})^{\varphi = 1} \text{ and } \mathcal{M}(T) \coloneqq (T\otimes_{\mathbf{Z}_p} \widehat{\mathcal{O}}_{\mathcal{E}}^{\mathrm{ur}})^{G_{\widetilde{L}_{\infty}}}.
\]  
The following equivalence is proved in \cite{kim-groupscheme-relative}.

\begin{prop}[cf.~{\cite[Prop.~7.7]{kim-groupscheme-relative}}] \label{prop:etale-phi-mod-Galois-rep}
The assignments $T(\cdot)$ and $\mathcal{M}(\cdot)$ are exact equivalences (quasi-inverse of each other) of $\otimes$-categories between $\mathrm{Mod}_{\mathcal{O}_{\mathcal{E}}}$ and $\mathrm{Rep}_{\mathbf{Z}_p}(G_{\widetilde{L}_{\infty}})$. Moreover, $T(\cdot)$ and $\mathcal{M}(\cdot)$ restrict to rank-preserving equivalence of categories between $\mathrm{Mod}_{\mathcal{O}_{\mathcal{E}}}^{\mathrm{free}}$ and $\mathrm{Rep}_{\mathbf{Z}_p}^{\mathrm{free}}(G_{\widetilde{L}_{\infty}})$ and length-preserving equivalence of categories between $\mathrm{Mod}_{\mathcal{O}_{\mathcal{E}}}^{\mathrm{tor}}$ and $\mathrm{Rep}_{\mathbf{Z}_p}^{\mathrm{tor}}(G_{\widetilde{L}_{\infty}})$. In both cases, $T(\cdot)$ and $\mathcal{M}(\cdot)$ commute with taking duals.   	
\end{prop}

For $\mathcal{M} \in \mathrm{Mod}_{\mathcal{O}_{\mathcal{E}}}^{\mathrm{free}}$ (resp. $\mathcal{M} \in \mathrm{Mod}_{\mathcal{O}_{\mathcal{E}}}^{\mathrm{tor}}$), we define $T^{\vee}(\mathcal{M}) \coloneqq \mathrm{Hom}_{\mathcal{O}_{\mathcal{E}}, \varphi}(\mathcal{M}, \widehat{\mathcal{O}}_{\mathcal{E}}^{\mathrm{ur}})$ (resp. $T^{\vee}(\mathcal{M}) \coloneqq \mathrm{Hom}_{\mathcal{O}_{\mathcal{E}}, \varphi}(\mathcal{M}, \widehat{\mathcal{O}}_{\mathcal{E}}^{\mathrm{ur}}\otimes_{\mathbf{Z}_p}\mathbf{Q}_p/ \mathbf{Z}_p)$). By the above Proposition, we have $T^{\vee}(\mathcal{M}) \cong (T(\mathcal{M}))^{\vee}$ in both cases. We similarly define $\mathcal{M}^{\vee}(T)$ for $T \in \mathrm{Rep}_{\mathbf{Z}_p}^{\mathrm{free}}(G_{\widetilde{L}_{\infty}})$ and $T \in \mathrm{Rep}_{\mathbf{Z}_p}^{\mathrm{tor}}(G_{\widetilde{L}_{\infty}})$. Suppose $0 \rightarrow \mathcal{M}_1 \rightarrow \mathcal{M}_2 \rightarrow \mathcal{M} \rightarrow 0$ is an exact sequence of \'etale $(\varphi, \mathcal{O}_{\mathcal{E}})$-modules with $\mathcal{M}_1, \mathcal{M}_2 \in \mathrm{Mod}_{\mathcal{O}_{\mathcal{E}}}^{\mathrm{free}}$ and $\mathcal{M} \in \mathrm{Mod}_{\mathcal{O}_{\mathcal{E}}}^{\mathrm{tor}}$. Then by Prop.~\ref{prop:etale-phi-mod-Galois-rep}, the corresponding sequence of representations
\[
0 \rightarrow T(\mathcal{M}_1) \rightarrow T(\mathcal{M}_2) \rightarrow T(\mathcal{M}) \rightarrow 0
\]
is exact. This induces the exact sequence of representations
\[
0 \rightarrow T^{\vee}(\mathcal{M}_2) \rightarrow T^{\vee}(\mathcal{M}_1) \rightarrow T^{\vee}(\mathcal{M}) \rightarrow 0,
\]
since $\mathrm{Hom}_{\mathbf{Z}_p}(T(\mathcal{M}), \mathbf{Z}_p) = 0 = \mathrm{Ext}^1_{\mathbf{Z}_p}(T(\mathcal{M}_2), \mathbf{Z}_p)$ and 
\[
T^{\vee}(\mathcal{M}) \cong (T(\mathcal{M}))^{\vee} = \mathrm{Hom}_{\mathbf{Z}_p}(T(\mathcal{M}), \mathbf{Q}_p/\mathbf{Z}_p) \cong \mathrm{Ext}^1_{\mathbf{Z}_p}(T(\mathcal{M}), \mathbf{Z}_p).
\]

For later use, we consider a natural base change as follows. Let $R_{0, g}$ be the $p$-completion of $\varinjlim_{\varphi} \mathcal{O}_{L_0}$, and let $k_g \coloneqq R_{0, g}/pR_{0, g}$. The Frobenius on $\mathcal{O}_{L_0}$ extends uniquely to $R_{0, g}$, and we have a $\varphi$-compatible isomorphism $R_{0, g} \cong W(k_g)$. Note that the map $\mathcal{O}_{L_0} \rightarrow W(k_g)$ is faithfully flat, since it sends $p$ to $p \in W(k_g)$ and both $\mathcal{O}_{L_0}$ and $W(k_g)$ are discrete valuation rings with $p$ being a uniformizer. Write $\mathfrak{S}_g \coloneqq W(k_g)[\![u]\!]$ and $\mathcal{O}_{\mathcal{E}, g} \coloneqq \mathfrak{S}_g[u^{-1}]^{\wedge}_p$ equipped with the induced Frobenius. Then the induced maps $\mathfrak{S} \rightarrow \mathfrak{S}_g$ and $\mathcal{O}_{\mathcal{E}} \rightarrow \mathcal{O}_{\mathcal{E}, g}$ are also $\varphi$-compatible and faithfully flat by the local criterion for flatness, and $\mathfrak{S}_g \cap \mathcal{O}_{\mathcal{E}} = \mathfrak{S}$ as subrings of $\displaystyle \mathcal{O}_{\mathcal{E}, g} = \{ \sum_{i \in \mathbf{Z}} a_iu^i ~|~ a_i \in W(k_g), \lim_{i \rightarrow -\infty}|a_i| \rightarrow 0 \}$ (where $|\cdot|$ denotes a $p$-adic norm on $W(k_g)$).

Let $K_g = W(k_g)\otimes_{W(k)} K$ and $K_{g, \infty} = \bigcup_{n \geq 0} K_g(\pi_n) \subset \overline{K_g}$. By \cite[Prop.~4.2.5 Pf.]{gao-integral-padic-hodge-imperfect}, under a suitable choice of embedding $\overline{L} \hookrightarrow \overline{K_g}$ extending $L \rightarrow K_g$ such that $t_{i, n} \in K_g$ for all $i$ and $n$, the map $G_{K_g} \coloneqq \mathrm{Gal}(\overline{K_g}/K_g) \rightarrow G_L$ induces $G_{K_{g, \infty}} \coloneqq \mathrm{Gal}(\overline{K_g}/K_{g, \infty}) \cong G_{\widetilde{L}_{\infty}}$. We fix such a choice. Note that the $p$-completions of $\overline{L}$ and $\overline{K_g}$ are identified, and so $W(\overline{L}^{\flat}) = W(\overline{K_g}^{\flat})$. Define $\widehat{\mathcal{O}}^{\mathrm{ur}}_{\mathcal{E}, g}$ similarly as above. For $\mathcal{M} \in \mathrm{Mod}_{\mathcal{O}_{\mathcal{E}}}$, $\mathcal{M}_g \coloneqq \mathcal{M}\otimes_{\mathcal{O}_{\mathcal{E}}} \mathcal{O}_{\mathcal{E}, g}$ equipped with the tensor-product Frobenius is an \'etale $(\varphi, \mathcal{O}_{\mathcal{E}, g})$-module. Furthermore, for $\mathcal{M}$ in either $\mathrm{Mod}_{\mathcal{O}_{\mathcal{E}}}^{\mathrm{free}}$ or $\mathrm{Mod}_{\mathcal{O}_{\mathcal{E}}}^{\mathrm{tor}}$, consider the natural $G_{K_{g, \infty}}$-equivariant map
\[
T(\mathcal{M}) \rightarrow T(\mathcal{M}_g) = (\mathcal{M}_g\otimes_{\mathcal{O}_{\mathcal{E}, g}} \widehat{\mathcal{O}}^{\mathrm{ur}}_{\mathcal{E}, g})^{\varphi = 1}.
\]
Note that we have natural isomorphisms
\begin{align*}
T(\mathcal{M})\otimes_{\mathbf{Z}_p} \widehat{\mathcal{O}}^{\mathrm{ur}}_{\mathcal{E}, g} &\cong (T(\mathcal{M})\otimes_{\mathbf{Z}_p} \widehat{\mathcal{O}}^{\mathrm{ur}}_{\mathcal{E}})\otimes_{\widehat{\mathcal{O}}^{\mathrm{ur}}_{\mathcal{E}}} \widehat{\mathcal{O}}^{\mathrm{ur}}_{\mathcal{E}, g} \\
	&\cong (\mathcal{M}\otimes_{\mathcal{O}_{\mathcal{E}}} \widehat{\mathcal{O}}^{\mathrm{ur}}_{\mathcal{E}})\otimes_{\widehat{\mathcal{O}}^{\mathrm{ur}}_{\mathcal{E}}} \widehat{\mathcal{O}}^{\mathrm{ur}}_{\mathcal{E}, g} \cong \mathcal{M}_g\otimes_{\mathcal{O}_{\mathcal{E}, g}} \widehat{\mathcal{O}}^{\mathrm{ur}}_{\mathcal{E}, g}\\
	&\cong T(\mathcal{M}_g)\otimes_{\mathbf{Z}_p} \widehat{\mathcal{O}}^{\mathrm{ur}}_{\mathcal{E}, g}.
\end{align*} 
Since $\mathbf{Z}_p \rightarrow \widehat{\mathcal{O}}^{\mathrm{ur}}_{\mathcal{E}, g}$ is faithfully flat, the above map $T(\mathcal{M}) \rightarrow T(\mathcal{M}_g)$ is an isomorphism.

\section{Kisin descent datum} \label{sec:prismatic-F-crystals}

We explain a result in \cite{du-liu-moon-shimizu-completed-prismatic-F-crystal-loc-system} classifying lattices of crystalline representations of $G_L$ via prismatic $F$-crystals on the absolute prismatic site of $\mathcal{O}_L$. We refer the readers to \cite[Def.~3.2]{bhatt-scholze-prismaticcohom-v3} for the terminologies on prisms. 

\begin{defn}[{\cite[Def.~2.3]{bhatt-scholze-prismaticFcrystal}}]  
The \textit{absolute prismatic site} $(\mathcal{O}_L)_{\Prism}$ of the $p$-adic formal scheme $\mathrm{Spf} \mathcal{O}_L$ consists of pairs $((A, I), \mathrm{Spf} A/I \rightarrow \mathrm{Spf} \mathcal{O}_L)$ where $(A, I)$ is a bounded prism and $\mathrm{Spf} A/I \rightarrow \mathrm{Spf} \mathcal{O}_L$ is a map of $p$-adic formal schemes. We often omit the structure map $\mathrm{Spf} A/I \rightarrow \mathrm{Spf} \mathcal{O}_L$ and simply write 	$(A, I) \in (\mathcal{O}_L)_{\Prism}$. The morphisms are the opposite of morphisms of bounded prisms over $\mathrm{Spf} \mathcal{O}_L$. We equip $(\mathcal{O}_L)_{\Prism}$ with the topology given by $(p, I)$-completely faithfully flat map of prisms $(A, I) \rightarrow (B, J)$ over $\mathrm{Spf} \mathcal{O}_L$.
\end{defn}

For any $(A, I) \in (\mathcal{O}_L)_{\Prism}$, $A$ is classically $(p, I)$-complete by \cite[Lem.~3.7]{bhatt-scholze-prismaticcohom-v3}. We have $(\mathfrak{S}, (E)) \in (\mathcal{O}_L)_{\Prism}$ with the structure map given by $\mathfrak{S}/(E) \cong \mathcal{O}_L$, called the \textit{Breuil--Kisin prism}. Denote by $(\mathfrak{S}^{(1)}, (E))$ (resp. $(\mathfrak{S}^{(2)}, (E))$) the self-product (resp. self-triple-product) of $(\mathfrak{S}, (E))$ in $(\mathcal{O}_L)_{\Prism}$ (see \cite[Ex.~3.4]{du-liu-moon-shimizu-completed-prismatic-F-crystal-loc-system}). The Breuil--Kisin prism covers the final object of $\mathrm{Shv}((\mathcal{O}_L)_{\Prism})$, so crystals on $(\mathcal{O}_L)_{\Prism}$ can be described by descent data involving the self-product and self-triple-product. Let $p_i\colon \mathfrak{S} \rightarrow \mathfrak{S}^{(1)}$, $i = 1, 2$, be the natural projection maps.

We first recall a result in \cite{bhatt-scholze-prismaticFcrystal} on $\mathbf{Z}_p$-representations of $G_L$.

\begin{defn}
Let $\mathrm{DD}_{\mathcal{O}_{\mathcal{E}}}$ denote the category consisting of tuples $(\mathcal{M}, \varphi_{\mathcal{M}}, f)$ where  	
\begin{itemize}
\item $(\mathcal{M}, \varphi_{\mathcal{M}}) \in \mathrm{Mod}_{\mathcal{O}_{\mathcal{E}}}^{\mathrm{free}}$;
\item $f\colon \mathfrak{S}^{(1)}[E^{-1}]^{\wedge}_p\otimes_{p_1, \mathfrak{S}[E^{-1}]^{\wedge}_p} \mathcal{M} \stackrel{\cong}{\rightarrow} \mathfrak{S}^{(1)}[E^{-1}]^{\wedge}_p\otimes_{p_2, \mathfrak{S}[E^{-1}]^{\wedge}_p} \mathcal{M}$ is an isomorphism of $\mathfrak{S}^{(1)}[E^{-1}]^{\wedge}_p$-modules compatible with Frobenii and satisfies the cocycle condition over $\mathfrak{S}^{(2)}[E^{-1}]^{\wedge}_p$.	
\end{itemize}
The morphisms of $\mathrm{DD}_{\mathcal{O}_{\mathcal{E}}}$ are $\mathcal{O}_{\mathcal{E}}$-module maps compatible with all structures.
\end{defn}

\begin{prop} \label{prop:equiv-GL-rep-laurent-F-crystal}
$\mathrm{DD}_{\mathcal{O}_{\mathcal{E}}}$ is naturally equivalent to $\mathrm{Rep}_{\mathbf{Z}_p}^{\mathrm{free}}(G_L)$. If $T(\mathcal{M}, \varphi_{\mathcal{M}}, f)$ denotes the object in $\mathrm{Rep}_{\mathbf{Z}_p}^{\mathrm{free}}(G_L)$ corresponding to $(\mathcal{M}, \varphi_{\mathcal{M}}, f) \in \mathrm{DD}_{\mathcal{O}_{\mathcal{E}}}$ under the equivalence, then $T(\mathcal{M}, \varphi_{\mathcal{M}}, f)|_{G_{\widetilde{L}_{\infty}}} = T(\mathcal{M})$ as representations of $G_{\widetilde{L}_{\infty}}$, where $T(\mathcal{M})$ is given as in Proposition~\ref{prop:etale-phi-mod-Galois-rep}. 	
\end{prop}

\begin{proof}
The equivalence of categories is a direct consequence of \cite[Cor.~3.8]{bhatt-scholze-prismaticFcrystal}, since the Breuil--Kisin prism covers the final object of $\mathrm{Shv}((\mathcal{O}_L)_{\Prism})$. Furthermore, by \cite[Thm.~3.2]{min-wang-rel-phi-gamma-prism-F-crys}, we have
\[
T(\mathcal{M}, \varphi_{\mathcal{M}}, f) = (\mathcal{M}\otimes_{\mathcal{O}_{\mathcal{E}}} W(\overline{L}^{\flat}))^{\varphi = 1}.
\]
The second statement then follows from the proof of \cite[Prop.~3.27 (ii)(b)]{du-liu-moon-shimizu-completed-prismatic-F-crystal-loc-system}: it is shown in \textit{loc. cit.} that the natural map
\[
T(\mathcal{M}) = (\mathcal{M}\otimes_{\mathcal{O}_{\mathcal{E}}}\widehat{\mathcal{O}}_{\mathcal{E}}^{\mathrm{ur}})^{\varphi = 1} \rightarrow (\mathcal{M}\otimes_{\mathcal{O}_{\mathcal{E}}} W(\overline{L}^{\flat}))^{\varphi = 1}
\]
induced by the embedding $\widehat{\mathcal{O}}_{\mathcal{E}}^{\mathrm{ur}} \rightarrow W(\overline{L}^{\flat})$ is bijective.
\end{proof}

By \cite[Cor.~4.36]{du-liu-moon-shimizu-completed-prismatic-F-crystal-loc-system}, $\mathbf{Z}_p$-lattices of crystalline representations of $G_L$ are classified by prismatic $F$-crystals on $(\mathcal{O}_L)_{\Prism}$. We will use an equivalent description in terms of \textit{Kisin descent datum}. Fix a positive integer $r$.

\begin{defn}
Let $\mathrm{Mod}_{\mathfrak{S}}^r$ denote the category consisting of tuples $(\mathfrak{M}, \varphi_{\mathfrak{M}})$ where 
\begin{itemize}
\item $\mathfrak{M}$ is a finite free $\mathfrak{S}$-module;
\item $\varphi_{\mathfrak{M}}\colon \mathfrak{M} \rightarrow \mathfrak{M}$ is a $\varphi$-semi-linear endomorphism such that $(\mathfrak{M}, \varphi_{\mathfrak{M}})$ has $E$-height $\leq r$, i.e. the cokernel of $1\otimes\varphi_{\mathfrak{M}}\colon \varphi^*\mathfrak{M} = \mathfrak{S}\otimes_{\varphi, \mathfrak{S}} \mathfrak{M} \rightarrow \mathfrak{M}$ is killed by $E^r$.
\end{itemize}
The morphisms are $\mathfrak{S}$-module maps compatible with Frobenius.	
\end{defn}

\begin{defn}[cf.~{\cite[Def.~3.24]{du-liu-moon-shimizu-completed-prismatic-F-crystal-loc-system}}]
Let $\mathrm{DD}_{\mathfrak{S}, [0, r]}$ denote the category consisting of triples $(\mathfrak{M}, \varphi_{\mathfrak{M}}, f)$ (called \textit{Kisin descent datum}) where
\begin{itemize}
\item $(\mathfrak{M}, \varphi_{\mathfrak{M}}) \in \mathrm{Mod}_{\mathfrak{S}}^r$;
\item $f\colon \mathfrak{S}^{(1)}\otimes_{p_1, \mathfrak{S}} \mathfrak{M} \stackrel{\cong}{\rightarrow} \mathfrak{S}^{(1)}\otimes_{p_2, \mathfrak{S}} \mathfrak{M}$ is an isomorphism of $\mathfrak{S}^{(1)}$-modules compatible with Frobenii and satisfies the cocycle condition over $\mathfrak{S}^{(2)}$.	
\end{itemize}
The morphisms of $\mathrm{DD}_{\mathfrak{S}, [0, r]}$ are $\mathfrak{S}$-module maps compatible with all structures.	
\end{defn}

We have a natural functor $\mathrm{DD}_{\mathfrak{S}, [0, r]} \rightarrow \mathrm{DD}_{\mathcal{O}_{\mathcal{E}}}$ given by extending the scalars $\mathfrak{M} \mapsto \mathfrak{M}\otimes_{\mathfrak{S}}\mathcal{O}_{\mathcal{E}}$ (equipped with the tensor-product Frobenius $\varphi_{\mathfrak{M}}\otimes\varphi_{\mathcal{O}_{\mathcal{E}}}$). Composing with the equivalence $\mathrm{DD}_{\mathcal{O}_{\mathcal{E}}} \cong \mathrm{Rep}_{\mathbf{Z}_p}^{\mathrm{free}}(G_L)$ in Proposition~\ref{prop:equiv-GL-rep-laurent-F-crystal} and taking the dual, we obtain a contravariant functor $\mathrm{DD}_{\mathfrak{S}, [0, r]} \rightarrow \mathrm{Rep}_{\mathbf{Z}_p}^{\mathrm{free}}(G_L)$. The following result is proved in \cite{du-liu-moon-shimizu-completed-prismatic-F-crystal-loc-system}.

\begin{thm}[cf.~{\cite[Prop.~3.25, Cor.~4.36]{du-liu-moon-shimizu-completed-prismatic-F-crystal-loc-system}}] \label{thm:equivalence-DD-cryst-rep}
The above functor induces an anti-equivalence between $\mathrm{DD}_{\mathfrak{S}, [0, r]}$ and the category of $\mathbf{Z}_p$-lattices of crystalline representations of $G_L$ with Hodge--Tate weights in $[0, r]$.	
\end{thm}

\section{Proof of main theorem}

In this section, we prove Theorem~\ref{thm:intro}. As explained in Section~\ref{sec:purity}, by Theorem~\ref{thm:purity}, it suffices to show the corresponding statement for representations of $G_L$. Fix a positive integer $r$. Let $T \in \mathrm{Rep}_{\mathbf{Z}_p}^{\mathrm{free}}(G_L)$. Suppose for each $n \geq 1$, $T/p^n T$ is \textit{torsion crystalline with Hodge--Tate weights in} $[0, r]$, i.e. there exist $G_L$-stable $\mathbf{Z}_p$-lattices $T_2^{(n)} \subset T_1^{(n)}$ of a crystalline representation with Hodge--Tate weights in $[0, r]$ such that $T/p^nT \cong T_1^{(n)} / T_2^{(n)}$ as $\mathbf{Z}_p[G_L]$-modules. We need to show that $T[p^{-1}]$ is a crystalline representation of $G_L$ with Hodge--Tate weights in $[0, r]$. 

The proof will proceed as follows. By Theorem~\ref{thm:equivalence-DD-cryst-rep}, it suffices to construct the Kisin descent datum $(\mathfrak{M}, \varphi_{\mathfrak{M}}, f) \in \mathrm{DD}_{\mathfrak{S}, [0, r]}$ corresponding to $T$. We first apply the analogous construction as in \cite[\S 4.3]{liu-fontaineconjecture} to obtain the $\mathfrak{S}$-module $\mathfrak{M}$ with Frobenius. Then $(\mathfrak{M}, \varphi_{\mathfrak{M}})$ being an object in $\mathrm{Mod}_{\mathfrak{S}}^r$ will follow from that $\mathfrak{S} \rightarrow \mathfrak{S}_g$ is faithfully flat and the results in \cite[\S 4.3]{liu-fontaineconjecture} for the $\mathfrak{S}_g$-module $\mathfrak{M}\otimes_{\mathfrak{S}} \mathfrak{S}_g$. 

To obtain the descent datum $f$, we consider
\[
f_{\text{\'et}}\colon \mathfrak{S}^{(1)}[E^{-1}]^{\wedge}_p\otimes_{p_1, \mathfrak{S}} \mathfrak{M} \stackrel{\cong}{\rightarrow} \mathfrak{S}^{(1)}[E^{-1}]^{\wedge}_p\otimes_{p_2, \mathfrak{S}} \mathfrak{M}
\]
associated with $T^{\vee}$ via Proposition~\ref{prop:equiv-GL-rep-laurent-F-crystal}. We will deduce from \cite[Cor.~3.7]{bhatt-scholze-prismaticFcrystal} together with some ring theoretic facts for $\mathfrak{S}^{(1)}$ that $f_{\text{\'et}}$ refines to $f\colon \mathfrak{S}^{(1)}\otimes_{p_1, \mathfrak{S}} \mathfrak{M} \stackrel{\cong}{\rightarrow} \mathfrak{S}^{(1)}\otimes_{p_2, \mathfrak{S}} \mathfrak{M}$. 

\begin{rem}
The argument in the above paragraph can also be used in the classical case $L = K$ and gives an alternate proof of the main result of \cite{liu-fontaineconjecture} for crystalline representations. Note that in \textit{loc. cit.}, after constructing $(\mathfrak{M}, \varphi_{\mathfrak{M}})$, some further computations for Galois invariants are made in \cite[\S 6--\S 8]{liu-fontaineconjecture}	 to show that $T[p^{-1}]$ is crystalline.
\end{rem}
 
We now begin the proof with the following observation.
 
\begin{lem} \label{lem:functor-fully-faithful}
The functor $\mathrm{Mod}_{\mathfrak{S}}^r \rightarrow \mathrm{Mod}_{\mathcal{O}_{\mathcal{E}}}$ given by $\mathfrak{M} \mapsto \mathfrak{M}\otimes_{\mathfrak{S}} \mathcal{O}_{\mathcal{E}}$ is fully faithful.
\end{lem}

\begin{proof}
Let $\mathfrak{M}_1, \mathfrak{M}_2 \in \mathrm{Mod}_{\mathfrak{S}}^r$, and let $f\colon \mathfrak{M}_1\otimes_{\mathfrak{S}} \mathcal{O}_{\mathcal{E}} \rightarrow \mathfrak{M}_2\otimes_{\mathfrak{S}} \mathcal{O}_{\mathcal{E}}$	be a morphism in $\mathrm{Mod}_{\mathcal{O}_{\mathcal{E}}}$. It suffices to show that $f(\mathfrak{M}_1) \subset \mathfrak{M}_2$. By \cite[Prop.~2.1.12]{kisin-crystalline}, we have $f(\mathfrak{M}_1) \subset \mathfrak{M}_2\otimes_{\mathfrak{S}} \mathfrak{S}_g$ as submodules of $\mathfrak{M}_2\otimes_{\mathfrak{S}} \mathcal{O}_{\mathcal{E}, g}$. Since $\mathfrak{M}_2$ is flat over $\mathfrak{S}$, we have
\[
f(\mathfrak{M}_1) \subset (\mathfrak{M}_2\otimes_{\mathfrak{S}} \mathcal{O}_{\mathcal{E}}) \cap (\mathfrak{M}_2\otimes_{\mathfrak{S}} \mathfrak{S}_g) = \mathfrak{M}_2\otimes_{\mathfrak{S}} (\mathcal{O}_{\mathcal{E}} \cap \mathfrak{S}_g) = \mathfrak{M}_2.
\] 
\end{proof}

Let $\mathcal{M} \coloneqq \mathcal{M}^{\vee}(T)$ be the \'etale $(\varphi, \mathcal{O}_{\mathcal{E}})$-module associated with $T$. Similarly, let $\mathcal{M}_i^{(n)} \coloneqq \mathcal{M}^{\vee}(T_i^{(n)})$ for $i = 1, 2$. The short exact sequence 
\[
0 \rightarrow T_2^{(n)} \rightarrow T_1^{(n)} \rightarrow T/p^n T \rightarrow 0
\]
induces a short exact sequence
\[
0 \rightarrow \mathcal{M}_1^{(n)} \rightarrow \mathcal{M}_2^{(n)} \stackrel{g_n}{\rightarrow} \mathcal{M}/p^n \mathcal{M} \rightarrow 0.
\]

For $i = 1, 2$, let $\mathfrak{M}_i^{(n)} \in \mathrm{DD}_{\mathfrak{S}, [0, r]}$ be the Kisin module of $E$-height $\leq r$ equipped with descent datum which corresponds to $T_i^{(n)}$ by Theorem~\ref{thm:equivalence-DD-cryst-rep}. In particular, $\mathfrak{M}_i^{(n)}\otimes_{\mathfrak{S}} \mathcal{O}_{\mathcal{E}} \cong \mathcal{M}_i^{(n)}$ as \'etale $(\varphi, \mathcal{O}_{\mathcal{E}})$-modules. By Lemma~\ref{lem:functor-fully-faithful}, we have induced maps
\[
\mathfrak{M}_1^{(n)} \hookrightarrow \mathfrak{M}_2^{(n)} \stackrel{g_n}{\rightarrow} \mathcal{M}/p^n\mathcal{M}.
\]
where the first map is injective. Write $\mathfrak{L} \coloneqq \ker(g_n\colon \mathfrak{M}_2^{(n)} \rightarrow \mathcal{M}/p^n\mathcal{M})$ and $\mathfrak{N} \coloneqq \mathrm{Im}(g_n) \subset \mathcal{M}/p^n\mathcal{M}$. We claim that the natural map $\mathfrak{M}_1^{(n)} \rightarrow \mathfrak{L}$ is an isomorphism, i.e. $\mathfrak{M}_1^{(n)} = \mathfrak{L}$ as $\mathfrak{S}$-submodules of $\mathcal{M}_1^{(n)}$. Since $\mathfrak{S} \rightarrow \mathfrak{S}_g$ is faithfully flat, it suffices to show that the induced map $\mathfrak{M}_1^{(n)}\otimes_{\mathfrak{S}} \mathfrak{S}_g \rightarrow \mathfrak{L}\otimes_{\mathfrak{S}} \mathfrak{S}_g$ is an isomorphism. Note first that $\mathfrak{M}_1^{(n)}\otimes_{\mathfrak{S}} \mathfrak{S}_g$ as a $\mathfrak{S}_g$-module with the induced Frobenius lies in $\mathrm{Mod}_{\mathfrak{S}_g}^r$. Furthermore, the cokernel of $1\otimes\varphi\colon \mathfrak{S}_g\otimes_{\varphi, \mathfrak{S}_g}(\mathfrak{N}\otimes_{\mathfrak{S}}\mathfrak{S}_g) \rightarrow \mathfrak{N}\otimes_{\mathfrak{S}}\mathfrak{S}_g$ is killed by $E^r$ by \cite[Prop.~1.3.5]{fontaine-p-adic-rep-I}. So we also have $\mathfrak{L}\otimes_{\mathfrak{S}} \mathfrak{S}_g \in \mathrm{Mod}_{\mathfrak{S}_g}^r$ by \cite[Cor.~2.3.8]{liu-fontaineconjecture}. On the other hand, since the maps $\mathcal{O}_{\mathcal{E}} \rightarrow \mathcal{O}_{\mathcal{E}, g}$ and $\mathfrak{S}_g \rightarrow \mathcal{O}_{\mathcal{E}, g}$ are flat, we have 
\[
(\mathfrak{M}_1^{(n)}\otimes_{\mathfrak{S}} \mathfrak{S}_g)\otimes_{\mathfrak{S}_g}\mathcal{O}_{\mathcal{E}, g} \cong \mathcal{M}_1^{(n)}\otimes_{\mathcal{O}_{\mathcal{E}}} \mathcal{O}_{\mathcal{E}, g} \cong (\mathfrak{L}\otimes_{\mathfrak{S}} \mathfrak{S}_g)\otimes_{\mathfrak{S}_g}\mathcal{O}_{\mathcal{E}, g}.
\]   
Thus, by \cite[Prop.~2.1.2]{kisin-crystalline}, we deduce $\mathfrak{M}_1^{(n)}\otimes_{\mathfrak{S}} \mathfrak{S}_g \cong \mathfrak{L}\otimes_{\mathfrak{S}} \mathfrak{S}_g$. 

We now have $\mathfrak{M}_2^{(n)}/\mathfrak{M}_1^{(n)} = g_n(\mathfrak{M}_2^{(n)}) \subset \mathcal{M}/p^n\mathcal{M}$. Denote $\mathfrak{M}_{(n)} \coloneqq \mathfrak{M}_2^{(n)}/\mathfrak{M}_1^{(n)}$. Note that $\mathfrak{M}_{(n)}\otimes_{\mathfrak{S}}\mathcal{O}_{\mathcal{E}} = \mathcal{M}/p^n\mathcal{M}$.

Let $e = [K : K_0]$ be the ramification index, and let $\mathfrak{c} \coloneqq \lfloor \frac{er}{p-1} \rfloor(2^{2r}(r\lfloor \frac{er}{p-1} \rfloor+1)-1)+1$ if $er \geq p-1$ and $\mathfrak{c} \coloneqq 0$ if $er < p-1$. For each $n \geq 1$, consider 
\begin{equation} \label{eq:M_n}
\mathfrak{M}_n \coloneqq \ker(p^{2\mathfrak{c}}\mathfrak{M}_{(n+3\mathfrak{c})} \stackrel{\times p^n}{\longrightarrow} p^{n+2\mathfrak{c}}\mathfrak{M}_{(n+3\mathfrak{c})}).
\end{equation}
Note that $\mathfrak{M}_i^{(n)}\otimes_{\mathfrak{S}} \mathfrak{S}_g \in \mathrm{Mod}_{\mathfrak{S}_g}^r$ for $i = 1, 2$. Since $\mathfrak{S} \rightarrow \mathfrak{S}_g$ is faithfully flat, we have $\mathfrak{M}_{(n)}\otimes_{\mathfrak{S}}\mathfrak{S}_g = (\mathfrak{M}_2^{(n)}\otimes_{\mathfrak{S}}\mathfrak{S}_g)/(\mathfrak{M}_1^{(n)}\otimes_{\mathfrak{S}}\mathfrak{S}_g)$ and 
\[
\mathfrak{M}_n\otimes_{\mathfrak{S}} \mathfrak{S}_g \cong \ker(p^{2\mathfrak{c}}(\mathfrak{M}_{(n+3\mathfrak{c})}\otimes_{\mathfrak{S}}\mathfrak{S}_g) \stackrel{\times p^n}{\longrightarrow} p^{n+2\mathfrak{c}}(\mathfrak{M}_{(n+3\mathfrak{c})}\otimes_{\mathfrak{S}}\mathfrak{S}_g)).
\]
By \cite[Lem.~4.4.1 Pf.]{liu-fontaineconjecture}, $\mathfrak{M}_n\otimes_{\mathfrak{S}} \mathfrak{S}_g$ is finite free over $\mathfrak{S}_g/p^n \mathfrak{S}_g$ and $p(\mathfrak{M}_{n+1}\otimes_{\mathfrak{S}} \mathfrak{S}_g) \cong \mathfrak{M}_n\otimes_{\mathfrak{S}} \mathfrak{S}_g$ compatibly with $\varphi$. Thus, since $\mathfrak{S}/p^n \mathfrak{S} \rightarrow \mathfrak{S}_g/p^n \mathfrak{S}_g$ is faithfully flat and $\mathfrak{S}/p^n \mathfrak{S}$ is local, we have that $\mathfrak{M}_n$ is finite free over $\mathfrak{S}/p^n\mathfrak{S}$ and $p\mathfrak{M}_{n+1} \cong \mathfrak{M}_n$. Furthermore, we deduce from the construction that $\mathfrak{M}_n\otimes_{\mathfrak{S}}\mathcal{O}_{\mathcal{E}} \cong \mathcal{M}/p^n\mathcal{M}$, and $\mathrm{coker}(1\otimes\varphi_{\mathfrak{M}_n}\colon \varphi^*\mathfrak{M}_n \rightarrow \mathfrak{M}_n)$ is killed by $E^r$. 

Let 
\[
\mathfrak{M} \coloneqq \varprojlim_n \mathfrak{M}_n
\] 
equipped with the induced Frobenius. Then $\mathfrak{M} \in \mathrm{Mod}_{\mathfrak{S}}^r$ and $\mathfrak{M}\otimes_{\mathfrak{S}} \mathcal{O}_{\mathcal{E}} \cong \mathcal{M}$ as \'etale $(\varphi, \mathcal{O}_{\mathcal{E}})$-modules. Furthermore, since $\mathfrak{M}_n$ is given by (\ref{eq:M_n}), we can apply the same construction as in \cite[Lem.~4.4.1 Pf.]{liu-fontaineconjecture} to obtain $G_L$-stable $\mathbf{Z}_p$-lattices $L_1^{(n)}, L_2^{(n)} \subset T_1^{(n)}[p^{-1}] = T_2^{(n)}[p^{-1}]$ sitting in an exact sequence of $\mathbf{Z}_p[G_L]$-modules
\[
0 \rightarrow L_2^{(n)} \rightarrow L_1^{(n)} \rightarrow T/p^n T \rightarrow 0
\]
such that by letting $\mathfrak{L}_i^{(n)} \in \mathrm{DD}_{\mathfrak{S}, [0, r]}$ be the Kisin module with descent datum corresponding to $L_i^{(n)}$ via Theorem~\ref{thm:equivalence-DD-cryst-rep}, we have the induced exact sequence
\[
0 \rightarrow \mathfrak{L}_1^{(n)} \rightarrow \mathfrak{L}_2^{(n)} \rightarrow \mathfrak{M}_n \rightarrow 0
\]
compatible with $\varphi$. 

We now construct a descent datum over $\mathfrak{S}^{(2)}$ for $\mathfrak{M}$. By Proposition~\ref{prop:equiv-GL-rep-laurent-F-crystal}, we have a natural $\mathfrak{S}^{(1)}$-linear isomorphism
\[
f_{\text{\'et}}\colon \mathfrak{S}^{(1)}[E^{-1}]^{\wedge}_p\otimes_{p_1, \mathfrak{S}} \mathfrak{M} \stackrel{\cong}{\rightarrow} \mathfrak{S}^{(1)}[E^{-1}]^{\wedge}_p\otimes_{p_2, \mathfrak{S}} \mathfrak{M}
\] 
associated with $T^{\vee}$, which is compatible with $\varphi$ and satisfies the cocycle condition over $\mathfrak{S}^{(2)}[E^{-1}]^{\wedge}_p$. The following is proved in \cite{du-liu-moon-shimizu-completed-prismatic-F-crystal-loc-system}.

\begin{lem}[cf.~{\cite[Cor.~3.6]{du-liu-moon-shimizu-completed-prismatic-F-crystal-loc-system}}] \label{lem:S2}
$\mathfrak{S}^{(1)}$ is $E$-torsion free and $p$-torsion free. Furthermore, $\mathfrak{S}^{(1)} = \mathfrak{S}^{(1)}[E^{-1}] \cap \mathfrak{S}^{(1)}[p^{-1}]$, and $\mathfrak{S}^{(1)}[E^{-1}]$ is $p$-adically separated. 	
\end{lem}

Since $\mathfrak{S}^{(1)} \rightarrow \mathfrak{S}^{(1)}[E^{-1}]^{\wedge}_p$ is injective by the above Lemma and $\mathfrak{M}$ is flat over $\mathfrak{S}$, the map $\mathfrak{S}^{(1)}\otimes_{p_i, \mathfrak{S}} \mathfrak{M} \rightarrow \mathfrak{S}^{(1)}[E^{-1}]^{\wedge}_p\otimes_{p_i, \mathfrak{S}} \mathfrak{M}$ for each $i = 1, 2$ is injective. 

\begin{prop} \label{prop:fet}
We have
\[
f_{\text{\'et}}(\mathfrak{S}^{(1)}\otimes_{p_1, \mathfrak{S}} \mathfrak{M}) \subset \mathfrak{S}^{(1)}\otimes_{p_2, \mathfrak{S}} \mathfrak{M}.
\]	
\end{prop}

\begin{proof}
From $\mathfrak{S}^{(1)} = \mathfrak{S}^{(1)}[E^{-1}] \cap \mathfrak{S}^{(1)}[p^{-1}]$ in Lemma~\ref{lem:S2}, it follows that the map 
\[
\mathfrak{S}^{(1)}/p^n\mathfrak{S}^{(1)} \rightarrow \mathfrak{S}^{(1)}[E^{-1}]^{\wedge}_p/p^n\mathfrak{S}^{(1)}[E^{-1}]^{\wedge}_p \cong \mathfrak{S}^{(1)}[E^{-1}]/p^n\mathfrak{S}^{(1)}[E^{-1}]
\]
is injective. So $\mathfrak{S}^{(1)} \cap p^n\mathfrak{S}^{(1)}[E^{-1}]^{\wedge}_p = p^n\mathfrak{S}^{(1)} \subset \mathfrak{S}^{(1)}[E^{-1}]^{\wedge}_p$. Since $\mathfrak{S}^{(1)}$ is classically $p$-complete, we have by \cite[Tag~031A]{stacks-project} that
\[
\bigcap_{n=1}^{\infty} (\mathfrak{S}^{(1)}+p^n\mathfrak{S}^{(1)}[E^{-1}]^{\wedge}_p) = \mathfrak{S}^{(1)}.
\]
Thus, it suffices to show
\[
f_{\text{\'et}}(\mathfrak{S}^{(1)}\otimes_{p_1, \mathfrak{S}} \mathfrak{M}) \subset \bigcap_{n=1}^{\infty} (\mathfrak{S}^{(1)}\otimes_{p_2, \mathfrak{S}} \mathfrak{M}+p^n\mathfrak{S}^{(1)}[E^{-1}]^{\wedge}_p\otimes_{p_2, \mathfrak{S}} \mathfrak{M})
\]
since $\mathfrak{M}$ is finite free over $\mathfrak{S}$. 

For each $n$, let 
\[
f_n\colon \mathfrak{S}^{(1)}\otimes_{p_1, \mathfrak{S}} \mathfrak{L}_{2}^{(n)} \stackrel{\cong}{\rightarrow} \mathfrak{S}^{(1)}\otimes_{p_2, \mathfrak{S}} \mathfrak{L}_{2}^{(n)}
\]
be the descent datum associated with $L_2^{(n)}$ via Theorem~\ref{thm:equivalence-DD-cryst-rep}. We deduce from \cite[Cor.~3.7]{bhatt-scholze-prismaticFcrystal} that the following diagram commutes:
\[
\xymatrix{
\mathfrak{S}^{(1)}[E^{-1}]^{\wedge}_p\otimes_{p_1, \mathfrak{S}} \mathfrak{L}_2^{(n)} \ar[r]^{f_n} \ar[d] & \mathfrak{S}^{(1)}[E^{-1}]^{\wedge}_p\otimes_{p_2, \mathfrak{S}} \mathfrak{L}_2^{(n)} \ar[d] \\
\mathfrak{S}^{(1)}[E^{-1}]^{\wedge}_p\otimes_{p_1, \mathfrak{S}} \mathfrak{M}/p^n\mathfrak{M} \ar[r]^{f_{\text{\'et}}} & \mathfrak{S}^{(1)}[E^{-1}]^{\wedge}_p\otimes_{p_2, \mathfrak{S}} \mathfrak{M}/p^n\mathfrak{M}
}
\]
Here, the vertical maps are surjections induced by $\mathfrak{L}_2^{(n)} \twoheadrightarrow \mathfrak{M}/p^n \mathfrak{M}$, and bottom horizontal map is given by $f_{\text{\'et}} \mod p^n$. Since $f_n(\mathfrak{S}^{(1)}\otimes_{p_1, \mathfrak{S}} \mathfrak{L}_2^{(n)}) \subset \mathfrak{S}^{(1)}\otimes_{p_2, \mathfrak{S}} \mathfrak{L}_2^{(n)}$ and the induced map $\mathfrak{S}^{(1)}\otimes_{p_1, \mathfrak{S}} \mathfrak{L}_2^{(n)} \rightarrow \mathfrak{S}^{(1)}\otimes_{p_1, \mathfrak{S}} \mathfrak{M}/p^n\mathfrak{M}$ is surjective, we have
\[
f_{\text{\'et}}(\mathfrak{S}^{(1)}\otimes_{p_1, \mathfrak{S}} \mathfrak{M}/p^n\mathfrak{M}) \subset \mathfrak{S}^{(1)}\otimes_{p_2, \mathfrak{S}} \mathfrak{M}/p^n\mathfrak{M}.
\]
This proves the assertion.  
\end{proof}

Applying Proposition~\ref{prop:fet} for $f_{\text{\'et}}^{-1}$, we obtain a descent datum 
\[
f\colon \mathfrak{S}^{(1)}\otimes_{p_1, \mathfrak{S}} \mathfrak{M} \stackrel{\cong}{\rightarrow} \mathfrak{S}^{(1)}\otimes_{p_2, \mathfrak{S}} \mathfrak{M}
\] 
which is compatible with $f_{\text{\'et}}$. Hence, $T$ is a $\mathbf{Z}_p$-lattice of crystalline representation of $G_L$ with Hodge--Tate weights in $[0, r]$ by Theorem~\ref{thm:equivalence-DD-cryst-rep}. 

This concludes the proof of the following: 

\begin{thm} \label{thm:main}
Fix a positive integer $r$, and let $T \in \mathrm{Rep}_{\mathbf{Z}_p}^{\mathrm{free}}(G_X)$. Suppose for each $n \geq 1$, there exist $G_X$-stable $\mathbf{Z}_p$-lattices $T_2^{(n)} \subset T_1^{(n)}$ of a crystalline local system on $X$ with Hodge--Tate weights in $[0, r]$ such that $T/p^nT \cong T_1^{(n)} / T_2^{(n)}$ as $\mathbf{Z}_p[G_X]$-modules. Then $T\otimes_{\mathbf{Z}_p} \mathbf{Q}_p$ is a crystalline local system on $X$ with Hodge--Tate weights in $[0, r]$.	
\end{thm}

\begin{rem}
The analogous statement for Barsotti--Tate representations does not hold in general by \cite[Prop.~5.8]{moon-relativeRaynaud}.
\end{rem}

\section{Crystalline deformation ring} \label{sec:cryst-deform-ring}

As an application of Theorem~\ref{thm:main}, we construct the quotient of the universal deformation ring parametrizing crystalline local systems on $X$ whose Hodge--Tate weights lie in a fixed interval.    

Let $E / \mathbf{Q}_p$ be a finite extension with residue field $\mathbf{F}$, and let $\mathcal{O}_E$ be its ring of integers. Let $\mathcal{C}$ denote the category of topological local $\mathcal{O}_E$-algebras $A$ such that:
\begin{itemize}
\item The map $\mathcal{O}_E \rightarrow A/\mathfrak{m}_A$ is surjective, where $\mathfrak{m}_A$ denotes the maximal ideal of $A$	.
\item The map from $A$ to the projective limit of its discrete artinian quotients is a topological isomorphism.
\end{itemize}

\noindent Note that the first condition implies $A/\mathfrak{m}_A = \mathbf{F}$. By the second condition, $A$ is complete and its topology is given by a collection of open ideals $\mathfrak{a}$ for which $A/\mathfrak{a}$ is artinian. Morphisms in $\mathcal{C}$ are continuous $\mathcal{O}_E$-algebra homomorphisms. 

Let $V_0$ be a continuous $\mathbf{F}$-representation of $G_X$ with rank $d$. Assume $V_0$ is absolutely irreducible. For $A \in \mathcal{C}$, a \textit{deformation} of $V_0$ in $A$ is an isomorphism class of continuous finite free $A$-representations $V$ of $G_X$ satisfying $\mathbf{F}\otimes_A V \cong V_0$ as $\mathbf{F}[G_X]$-modules. We denote by $\mathrm{Def}(V_0, A)$ the set of such deformations. A morphism $A \rightarrow A'$ in $\mathcal{C}$ induces a map $f_*\colon \mathrm{Def}(V_0, A) \rightarrow \mathrm{Def}(V_0, A')$ sending the class of a representation $V$ over $A$ to the class of $A'\otimes_{f, A} V$. The following is proved in \cite{smit-lenstra}.

\begin{thm}[cf.~{\cite[Thm.~2.3]{smit-lenstra}}]
There exists a universal deformation ring $A_{\mathrm{univ}} \in \mathcal{C}$ and a deformation $V_{A_{\mathrm{univ}}} \in \mathrm{Def}(V_0, A_{\mathrm{univ}})$ such that for all rings $A \in \mathcal{C}$, we have a bijection 
\begin{equation} \label{eq:univ-deform}
\mathrm{Hom}_{\mathcal{C}}(A_{\mathrm{univ}}, A) \stackrel{\cong}{\rightarrow} \mathrm{Def}(V_0, A)
\end{equation}
given by $f \mapsto f_*(V_{A_{\mathrm{univ}}})$. 
\end{thm}

Let $\mathcal{C}^0$ be the full subcategory of $\mathcal{C}$ consisting of artinian rings. Abusing the notation, we write $V \in \mathrm{Def}(V_0, A)$ for a continuous $A$-representation $V$ of $G_X$ to mean $\mathbf{F}\otimes_A V \cong V_0$. Fix an integer $r \geq 1$. Denote by $\mathrm{Rep}_{\mathbf{Q}_p}^{\mathrm{cris, [0, r]}}(G_X)$ the category of \'etale isogeny $\mathbf{Z}_p$-local systems on $X$ which are crystalline with Hodge--Tate weights in $[0, r]$. For $A \in \mathcal{C}^0$ and $V_A \in \mathrm{Def}(V_0, A)$, we say $V_A$ is \textit{crystalline with Hodge--Tate weights in} $[0, r]$ if there exist a finite flat $\mathcal{O}_E$-algebra $B$, a surjection $g\colon B \twoheadrightarrow A$ of $\mathcal{O}_E$-algebras, and a continuous finite free $B$-representation $V_B$ of $G_X$ such that $V_B \otimes_{\mathbf{Z}_p} \mathbf{Q}_p \in \mathrm{Rep}_{\mathbf{Q}_p}^{\mathrm{cris, [0, r]}}(G_X)$ and $A \otimes_{g, B} V_B \cong V_A$ as $A[G_X]$-modules. For $A \in \mathcal{C}$, denote by $S_{\mathrm{cris}, [0, r]}(A)$ the subset of $\mathrm{Def}(V_0, A)$ consisting of the isomorphism classes of representations $V_A$ such that $A/\mathfrak{a}\otimes_A V_A$ is crystalline with Hodge--Tate weights in $[0, r]$ for all open ideals $\mathfrak{a} \subsetneq A$.

\begin{prop} \label{prop:quot-deform}
For any $\mathcal{C}$-morphism $f\colon A \rightarrow A'$, we have $f_*(S_{\mathrm{cris}, [0, r]}(A)) \subset S_{\mathrm{cris}, [0, r]}(A')$. There exists a closed ideal $\mathfrak{a}_{\mathrm{cris}, [0, r]}$ of $A_{\mathrm{univ}}$ such that the map (\ref{eq:univ-deform}) induces a bijection $\mathrm{Hom}_{\mathcal{C}}(A_{\mathrm{univ}}/\mathfrak{a}_{\mathrm{cris}, [0, r]}, A) \stackrel{\cong}{\rightarrow} S_{\mathrm{cris}, [0, r]}(A)$.
\end{prop}

\begin{proof}
We need to check the conditions in \cite[\S 6]{smit-lenstra}. Let $f\colon A \hookrightarrow A'$ be an inclusion of artinian rings in $\mathcal{C}$, and let $V_A \in \mathrm{Def}(V_0, A)$. We first claim that $V_A \in S_{\mathrm{cris}, [0, r]}(A)$ if and only if $V_{A'} \coloneqq A'\otimes_{f, A} V_A \in S_{\mathrm{cris}, [0, r]}(A')$. Suppose $V_A \in S_{\mathrm{cris}, [0, r]}(A)$. Then there exist a finite flat $\mathcal{O}_E$-algebra $B$, a surjection $g\colon B \twoheadrightarrow A$, and a finite free $B$-representation $V_B$ of $G_X$ such that $V_B \otimes_{\mathbf{Z}_p} \mathbf{Q}_p \in \mathrm{Rep}_{\mathbf{Q}_p}^{\mathrm{cris, [0, r]}}(G_X)$ and $A \otimes_{g, B} V_B \cong V_A$. 

There exists a surjection $f'\colon A[x_1, \ldots, x_n] \twoheadrightarrow A'$ of $\mathcal{O}_E$-algebras extending $f$ such that $f'(x_i) \in \mathfrak{m}_{A'}$ for each $i$. Let $I_{m, A} \subset A[x_1, \ldots, x_n]$ denote the ideal generated by the $m$-th degree homogeneous polynomials with coefficients in $A$. Since $A'$ is artinian, $f'(I_{m, A}) = 0$ for a sufficiently large $m$, and $f'$ induces a surjection $A[x_1, \ldots, x_n]/I_{m, A} \twoheadrightarrow A'$ for such $m$. Thus, we have surjective homomorphisms of $\mathcal{O}_E$-algebras
\[
g'\colon B' \coloneqq B[x_1, \ldots, x_n]/ I_{m, B} \twoheadrightarrow A[x_1, \ldots, x_n]/I_{m, A} \twoheadrightarrow A'.
\]   
Note that $B'$ is a finite flat $\mathcal{O}_E$-algebra. Let $V_{B'} = B'\otimes_B V_B$. As explained in Section~\ref{sec:purity} using Theorem~\ref{thm:purity}, to show $V_{B'}\otimes_{\mathbf{Z}_p}\mathbf{Q}_p \in \mathrm{Rep}_{\mathbf{Q}_p}^{\mathrm{cris, [0, r]}}(G_X)$, it suffices to consider the case of representations of $G_L$, which follows from Lemma~\ref{lem:base-change-coeff-reps}. We have $A'\otimes_{g', B'}V_{B'} \cong V_{A'}$, so $V_{A'} \in S_{\mathrm{cris}, [0, r]}(A')$. 	

Conversely, suppose $V_{A'} \in S_{\mathrm{cris}, [0, r]}(A')$. There exist a finite flat $\mathcal{O}_E$-algebra $B'$, a surjection $g'\colon B' \twoheadrightarrow A'$, and a finite free $B'$-representation $V_{B'}$ of $G_X$ such that $V_{B'} \otimes_{\mathbf{Z}_p} \mathbf{Q}_p \in \mathrm{Rep}_{\mathbf{Q}_p}^{\mathrm{cris, [0, r]}}(G_X)$ and $A' \otimes_{g', B'} V_{B'} \cong V_{A'}$ as $A'$-representations. Let $B$ be the kernel of the composite of morphisms $B' \stackrel{g'}{\rightarrow} A' \rightarrow A'/f(A)$ of $\mathcal{O}_E$-modules. Note that $B$ has an $\mathcal{O}_E$-algebra structure induced from the inclusion $B \hookrightarrow B'$, since $f\colon A \rightarrow A'$ is an $\mathcal{O}_E$-algebra map. Furthermore, since $B \subset B'$ which is $p$-torsion free, $B$ is a finite flat $\mathcal{O}_E$-algebra, and we have the surjection $g\colon B \twoheadrightarrow A$ of $\mathcal{O}_E$-algebras induced from $g'$. Let $V_B$ be the kernel of the following composite
\[
V_{B'} \rightarrow A' \otimes_{g', B'} V_{B'} \stackrel{\cong}{\rightarrow} V_{A'} \rightarrow A'/f(A)\otimes_{A'}V_{A'}.
\]
Then by construction, $V_B$ is a continuous finite free $B$-representation of $G_X$ such that $B'\otimes_B V_B \cong V_{B'}$. Furthermore, we have $A\otimes_{g, B} V_B \cong V_A$ induced from $A' \otimes_{g', B'} V_{B'} \cong V_{A'}$, since $V_{A'} = A'\otimes_{f,A} V_A$. By reducing to the case of representations of $G_L$ and using Lemma~\ref{lem:exact-seq-reps} and Lemma~\ref{lem:base-change-coeff-reps}, we deduce that $V_B \otimes_{\mathbf{Z}_p}\mathbf{Q}_p \in \mathrm{Rep}_{\mathbf{Q}_p}^{\mathrm{cris, [0, r]}}(G_X)$. So $V_A \in S_{\mathrm{cris}, [0, r]}(A)$.

Now, for $A \in \mathcal{C}$ and $V_A \in \mathrm{Def}(V_0, A)$, suppose $\mathfrak{a}_1, \mathfrak{a}_2 \subsetneq A$ are open ideals such that $A/\mathfrak{a}_i\otimes_A V_A \in S_{\mathrm{cris}, [0, r]}(A/\mathfrak{a}_i)$ for $i = 1, 2$. We claim that $A/(\mathfrak{a}_1\cap\mathfrak{a}_2)\otimes_A V_A \in S_{\mathrm{cris}, [0, r]}(A/(\mathfrak{a}_1\cap\mathfrak{a}_2))$. There exist a finite flat $\mathcal{O}_E$-algebra $B_i$, a surjection $g_i: B_i \twoheadrightarrow A/\mathfrak{a}_i$, and a finite free $B_i$-representation $V_{B_i}$ of $G_X$ such that $V_{B_i}\otimes_{\mathbf{Z}_p}\mathbf{Q}_p \in \mathrm{Rep}_{\mathbf{Q}_p}^{\mathrm{cris, [0, r]}}(G_X)$ and $A/\mathfrak{a}_i\otimes_{g_i, B_i}V_{B_i} \cong A/\mathfrak{a}_i\otimes_A V_A$. Let $V_{B_1 \times B_2}$ be the $(B_1 \times B_2)$-representation corresponding to $V_{B_1}\oplus V_{B_2}$. Note that $V_{B_1 \times B_2}\otimes_{\mathbf{Z}_p}\mathbf{Q}_p \in \mathrm{Rep}_{\mathbf{Q}_p}^{\mathrm{cris, [0, r]}}(G_X)$. Consider the inclusion $A/(\mathfrak{a}_1\cap\mathfrak{a}_2) \subset A/\mathfrak{a}_1 \times A/\mathfrak{a}_2$. We will follow a similar argument as in the previous paragraph. Let $B$ be the kernel of the composite 
\[
B_1 \times B_2 \stackrel{g_1\times g_2}{\longrightarrow} A/\mathfrak{}a_1 \times A/\mathfrak{a}_2 \rightarrow (A/\mathfrak{a}_1 \times A/\mathfrak{a}_2)/(A/(\mathfrak{a}_1\cap\mathfrak{a}_2)) 
\]   
of $\mathcal{O}_E$-module maps. Similarly as in the previous paragraph, $B$ is a finite flat $\mathcal{O}_E$-algebra, and we have the surjection $g\colon B \rightarrow A/(\mathfrak{a}_1\cap\mathfrak{a}_2)$ induced from $g_1 \times g_2$. Let $V_B$ be the kernel of the composite of morphisms
\[
V_{B_1 \times B_2} \rightarrow (A/\mathfrak{a}_1 \times A/\mathfrak{a}_2)\otimes_{g_1\times g_2, B_1\times B_2}V_{B_1 \times B_2} \cong (A/\mathfrak{a}_1 \times A/\mathfrak{a}_2)\otimes_A V_A
\]
and
\[
(A/\mathfrak{a}_1 \times A/\mathfrak{a}_2)\otimes_A V_A \rightarrow (A/\mathfrak{a}_1 \times A/\mathfrak{a}_2)/(A/(\mathfrak{a}_1\cap\mathfrak{a}_2)) \otimes_A V_A.
\]
By construction, $V_B$ is a continuous finite free $B$-representation of $G_X$ such that $(B_1 \times B_2) \otimes_B V_B \cong V_{B_1\times B_2}$ and $A/(\mathfrak{a}_1\cap\mathfrak{a}_2)\otimes_{g, B}V_B \cong A/(\mathfrak{a}_1 \cap \mathfrak{a}_2) \otimes_A V_A$. Similarly as above, we obtain $V_B \otimes_{\mathbf{Z}_p}\mathbf{Q}_p \in \mathrm{Rep}_{\mathbf{Q}_p}^{\mathrm{cris, [0, r]}}(G_X)$ by reducing to the case of representations of $G_L$ and using Lemma~\ref{lem:exact-seq-reps} and \ref{lem:base-change-coeff-reps}. Thus, $A/(\mathfrak{a}_1 \cap \mathfrak{a}_2)\otimes_A V_A \in S_{\mathrm{cris}, [0, r]}(A/(\mathfrak{a}_1 \cap \mathfrak{a}_2))$.

The assertion then follows from \cite[Prop.~6.1]{smit-lenstra}.
\end{proof}

\begin{thm} \label{thm:cryst-deform-ring}
Let $A$ be a finite flat $\mathcal{O}_E$-algebra, and let $f\colon A_{\mathrm{univ}} \rightarrow A$ be a continuous $\mathcal{O}_E$-algebra homomorphism. The induced representation $A[p^{-1}]\otimes_{f, A_{\mathrm{univ}}} V_{A_{\mathrm{univ}}}$ is a crystalline local system on $X$ with Hodge--Tate weights in $[0, r]$ if and only if $f$ factors through the quotient $A_{\mathrm{univ}}/\mathfrak{a}_{\mathrm{cris}, [0, r]}$.	
\end{thm}

\begin{proof}
Write $A_1 = f(A_{\mathrm{univ}}) \subset A$. Note that $A_1$ is local and finite flat over $\mathcal{O}_E$. Equip $A_1$ with the $p$-adic topology, so that $A_1 \in \mathcal{C}$ and $f\colon A_{\mathrm{univ}} \rightarrow A_1$ is a morphism in $\mathcal{C}$. Denote $V_{A_1} = A_1\otimes_{f, A_{\mathrm{univ}}} V_{A_{\mathrm{univ}}}$ and $V_{A} = A \otimes_{f, A_{\mathrm{univ}}} V_{A_{\mathrm{univ}}}$. Note that $V_A = A\otimes_{A_1} V_{A_1}$.

First suppose $V_A[p^{-1}]$ is a crystalline local system on $X$ with Hodge--Tate weights in $[0, r]$. By reducing to the case of representations of $G_L$ and using Lemma~\ref{lem:exact-seq-reps} and \ref{lem:base-change-coeff-reps}, we deduce $V_{A_1}[p^{-1}] \in \mathrm{Rep}_{\mathbf{Q}_p}^{\mathrm{cris, [0, r]}}(G_X)$. So $f$ factors through $A_{\mathrm{univ}}/\mathfrak{a}_{\mathrm{cris}, [0, r]}$ by Proposition~\ref{prop:quot-deform}.

Conversely, suppose $f$ factors through $A_{\mathrm{univ}}/\mathfrak{a}_{\mathrm{cris}, [0, r]}$. By Proposition~\ref{prop:quot-deform}, we have $V_{A_1} \in S_{\mathrm{cris}, [0, r]}(A_1)$, so $A_1/p^n\otimes_{A_1} V_{A_1}$ is torsion crystalline with Hodge--Tate weights in $[0, r]$ for each $n \geq 1$. By Theorem~\ref{thm:main}, $V_{A_1}[p^{-1}] \in \mathrm{Rep}_{\mathbf{Q}_p}^{\mathrm{cris, [0, r]}}(G_X)$, and so $V_A[p^{-1}] \in \mathrm{Rep}_{\mathbf{Q}_p}^{\mathrm{cris, [0, r]}}(G_X)$.      
\end{proof}

\bibliographystyle{amsalpha}
\bibliography{library}
	
\end{document}